# Frobenius pull-back and stability of vector bundles in characteristic 2

Jiu-Kang Yu[1]

and

Eugene Z. Xia

National Center for Theoretical Sciences
Hsinchu, Taiwan

November 6, 2001

**1 Introduction.** Let $X$ be an irreducible smooth projective curve of genus $g$ over an algebraically closed field $k$ of characteristic $p > 0$, and $F : X \to X$ the absolute Frobenius morphism on $X$. It is known that pulling back a stable vector bundle on $X$ by $F$ may destroy stability. One may measure the failure of (semi-)stability by the Harder-Narasimhan polygons of vector bundles.

In more formal language, let $n \geq 2$ be an integer, $\mathfrak{M}$ the coarse moduli space of stable vector bundles of rank $n$ and a fixed degree on $X$. Applying a theorem of Shatz (an analogue of the Grothendieck specialization theorem for $F$-isocrystals) to the pull-back by $F$ of the universal bundle (assuming the existence) on $\mathfrak{M}$, we see that $\mathfrak{M}$ has a canonical stratification by Harder-Narasimhan polygons ([S]). This interesting extra structure on $\mathfrak{M}$ is a feature of characteristic $p$. However, very little is known about these strata. Scattered constructions of points outside of the largest (semi-stable) stratum can be found in [G], [RR], [R], and [JX].

This paper deals exclusively with $p = 2$ and $n = 2$. On any curve $X$ of genus $\geq 2$, we provide a complete classification of rank-2 semi-stable vector bundles $V$ with $F^*V$ not semi-stable. We also obtain fairly good information about the locus destablized by Frobenius in the moduli space, including the irreducibility and the dimension of each non-empty Harder-Narasimhan stratum. This shows that the bound in [Su, Theorem 3.1] is sharp. A very interesting consequence of our classification is that high unstability of $F^*V$ implies high stability of $V$.

We conclude this introduction by remarking that the problem studied here can be cast in the generality of principal $G$-bundles over $X$, where $G$ is a connected reductive group over $k$. More precisely, consider the pull-back by $F$ of the universal object on the moduli stack of semi-stable principal $G$-bundles on $X$. Atiyah-Bott's generalization of the Harder-Narasimhan filtration should then give a canonical stratification of the moduli stack ([AB], see also [C]). In this context, our paper treats the case of $p = 2$, $G = \mathrm{GL}_2$.

We benefited from discussions with C.-L. Chai and K. Joshi. Both authors thank the hospitality of the National Center for Theoretical Sciences, Hsinchu, Taiwan, where most of this work was done during our visits in 2000-2001.

---

[1]partially supported a Sloan Foundation Fellowship and by grant DMS 0100678 from the National Science Foundation

**2  A measure of stability.**  In this paper, "a vector bundle" always means "a vector bundle over $X$". Following [LN], for a rank-2 vector bundle $V$, we put

$$s(V) = \deg(V) - 2\max\{\deg(L) : L \hookrightarrow V\},$$

where the maximum is taken over all rank-1 sub-module of $V$. By definition, $s(V) > 0$ (resp. $s(V) \geq 0$) if and only if $V$ is stable (resp. semi-stable). When $s(V) \leq 0$, the information of $(s(V), \deg(V))$ is the same as that of the Harder-Narasimhan polygon of $V$. Therefore, one may regard $s$ as a measure of stability extrapolating the Harder-Narasimhan polygons, though it is only for the rank-2 case (for possible variants for the higher rank case, see [BL]; for general reductive group, see [HN]).

**3  Raynaud's distinguished theta characteristic.**  From now on, $p = 2$ and $g \geq 2$. Following Raynaud [R, §4], we have a distinguished line bundle $B$ on $X$ defined by the exact sequence

$$0 \to \mathcal{O}_X \to F_*\mathcal{O}_X \to B \to 0.$$

Raynaud shows that $B^2 \simeq \Omega_X$, i.e. $B$ is a theta characteristic.

**Proposition.** *Let $\xi$ be a line bundle on $X$ and put $V = F_*(\xi \otimes B^{-1})$. Then $V$ is a vector bundle of rank 2 such that $\det(V) = \xi$, $V$ is stable, and $F^*V$ is not semi-stable. In fact,*

$$s(V) \geq g - 1 \text{ and } s(F^*V) = -(2g - 2).$$

*If $M$ is the sub-bundle of $F^*V$ of rank 1 such that $\deg M > \deg(F^*V)/2$, then $M \simeq \xi B$.*

PROOF. Write $L = \xi \otimes B^{-1}$. On an affine open set $U$ on which $F_*\mathcal{O}_X$, $B$, $L$ are trivial, choose a section $s \in (F_*\mathcal{O}_X)(U)$ such that the image of $s$ generates $B = F_*(\mathcal{O}_X)/\mathcal{O}_X$ on $U$, and a section $t \in L(U)$ generating $L|U$. Then $\{t, st\}$ generates $(F_*L)|U$ and $t \wedge (st) \mapsto t \otimes s$ is an isomorphism $\det(F_*L)|U \to (L \otimes B)|U$. One can check that this isomorphism is independent of the choices of $s, t$; hence, we obtain an isomorphism $\det(F_*L) \to L \otimes B$ by gluing these isomorphisms over various $U$'s.

Write $L = \xi \otimes B^{-1}$ and $d = \deg L$. Notice that $\deg V = d + g - 1$. Suppose that $M \hookrightarrow V$ is a sub-bundle of rank 1. By adjunction, there is a non-zero morphism $F^*M \to L$. Therefore, $\deg(F^*M) \leq \deg(L)$. Thus $\deg M \leq d/2 < (d + g - 1)/2 = \deg(V)/2$. Therefore, $V$ is stable and $s(V) \geq g - 1$.

Now consider the identity morphism $F_*L \to F_*L$. By adjunction, this gives a non-zero morhism $F^*V \to L$, which is surjective by a local calculation. The kernel of this morphism is a line bundle of degree $2(d + g - 1) - d = d + 2g - 2 > d + g - 1 = \deg(F^*V)/2$. So $F^*V$ is not semi-stable and $s(F^*V) = -(2g - 2)$. ∎

**Remark.** Let $V = F_*(\xi \otimes B^{-1})$. The extension

(*) $$0 \to \xi \otimes B \to F^*V \to \xi \otimes B^{-1} \to 0$$

defines a class in $\mathrm{Ext}^1(\xi \otimes B^{-1}, \xi \otimes B) \simeq H^1(X, B^2) \simeq k$. This class is trivial precisely when $\deg(\xi \otimes B^{-1})$ is even.

PROOF. Suppose that $\deg(\xi \otimes B^{-1})$ is even. Then we can write $L = \xi \otimes B^{-1} = M^2$. By [JX, §2], there is an exact sequence $0 \to M \to V \to M \otimes B \to 0$. Pulling back by $F$, we get $0 \to L \to F^*V \to L \otimes B^2 \to 0$. This shows that $(*)$ is split.

Suppose that $L = \xi \otimes B^{-1}$ has odd degree $2n+1$. By a theorem of Nagata ([LN], Cf. 8) and the above proposition, there is an exact sequence $0 \to M_1 \to V \to M_2 \to 0$, where $M_1, M_2$ are line bundles with degrees $n$ and $n+g$ respectively. From the exact sequence $0 \to M_1^2 \to F^*V \to M_2^2 \to 0$, we deduce that $\dim \mathrm{Hom}(L, F^*V) \leq \dim \mathrm{Hom}(L, M_1^2) + \dim \mathrm{Hom}(L, M_2^2) = 0 + g = g$ by the Riemann-Roch formula. Since $\mathrm{Hom}(L, \xi \otimes B) = H^0(X, B^2)$ has dimension $g$, any morphism $L \to F^*V$ factors through the sub-module $\xi \otimes B$ in $(*)$. Therefore, $(*)$ is not split. ∎

**4 The basic construction.** Henceforth, fix an integer $d$. For an injection $V' \hookrightarrow V''$ of vector bundles of the same rank, define the *co-length* $l$ of $V'$ in $V''$ to be the length of the torsion $\mathcal{O}_X$-module $V''/V'$. Clearly, $s(V') \geq s(V'') - l$.

We now give a basic construction of stable vector bundles $V$ of rank 2 with $F^*V$ not semi-stable. Let $l \leq g - 2$ be a non-negative integer, $L$ a line bundle of degree $d - 1 - (g - 2 - l)$, and $V$ a sub-module of $F_*L$ of co-length $l$, then $\deg V = d$ and $s(V) \geq (g-1) - l > 0$ by by Proposition 3. Therefore, $V$ is stable.

On the other hand, by adjunction, there is a morphism $F^*V \to L$, and the kernel is a line bundle of degree $\geq d + 1 + (g - 2 - l) > d = \deg(F^*V)/2$. Therefore, $F^*V$ is not semi-stable.

**5 Exhaustion.** Suppose that $V$ is semi-stable of rank 2 and $F^*V$ is not semi-stable.

Let $\xi = \det(V)$ and $d = \deg \xi = \deg V$. Since $F^*V$ is not semi-stable and of degree $2d$, there are line bundles $L, L'$ and an exact sequence $0 \to L' \to F^*V \to L \to 0$ with $\deg L' \geq d+1$, $\deg L \leq d-1$. By adjunction, this provides a non-zero morphism $V \to F_*L$. If the image is a line bundle $M$, we have $\deg M \geq d/2$ by semi-stability of $V$, and $\deg M \leq (d - 1 + g - 1)/2 - (g-1)/2 = (d-1)/2$ by Proposition 3. This is a contradiction.

Thus the image is of rank 2. Since $\deg V = d$ and $\deg(F_*L) \leq d + (g-2)$, $V$ is a sub-module of $F_*L$ of co-length $l \leq g - 2$, and $\deg L = d - 1 - (g - 2 - l)$.

Thus the basic construction yields all semi-stable vector bundles $V$ of rank 2, with $F^*V$ not semi-stable.

**5.1 Corollary.** *If $V$ is semi-stable of rank 2 with $F^*V$ not semi-stable, then $V$ is actually stable.* ∎

**5.2 Corollary.** *The basic construction with $l = g - 2$ already yields all semi-stable vector bundles $V$ of rank 2, with $F^*V$ not semi-stable.*

PROOF. In fact, if $l < l' \leq g - 2$ and $L' = L \otimes \mathcal{O}(D)$ for some effective divisor $D$ of degree $l' - l$ on $X$, then $V \hookrightarrow F_*L \hookrightarrow F_*L'$. Hence $V$ is also a sub-module of $F_*L'$ of co-length $l'$. Thus $V$ arises from the basic construction with $(l', L')$ playing the role of $(l, L)$. ∎

**6 Classification.** Let $L$ be a line bundle and let $Q = Q_l = Q_{l,L} = \text{Quot}_l(F_*L/X/k)$ be the scheme classifying sub-modules of $F_*L$ of co-length $l$ ([FGA, 3.2]). Let

$$\mathcal{V} \hookrightarrow \mathcal{O}_Q \boxtimes F_*L = (\text{id} \times F)_*(\mathcal{O}_Q \boxtimes L)$$

(sheaves on $Q \times X$) be the universal object on $Q$. By adjunction, we have a morphism $(1 \times F)^*\mathcal{V} \to \mathcal{O}_Q \boxtimes L$. Let $\mathcal{F}$ be the cokernel. Then $\text{pr}_*\mathcal{F}$ is a coherent sheaf on $Q$, where $\text{pr} : Q \times X \to Q$ is the projection ([H, II.5.20]). By [H, III.12.7.2], the subset

$$\{q \in Q : \dim_{\kappa(q)}((\text{pr}_*\mathcal{F}) \otimes \kappa(q) > 0\}$$

is closed. Its complement is an open sub-scheme, denoted by $Q^* = Q_l^* = Q_{l,L}^*$, of $Q$. Then $Q^*$ parametrizes those $V$'s with surjective $F^*V \to L$.

Let $\overline{\mathcal{M}}$ be the coarse moduli space of rank-2 semi-stable vector bundles of degree $d$ on $X$. Let $\mathcal{M}$ be the open sub-scheme parametrizing stable vector bundles and and $\mathcal{M}_1(k) \subset \overline{\mathcal{M}}(k)$ the subset of those $V$'s such that $F^*V$ is not semi-stable. By Corollary 5.1, $\mathcal{M}_1(k) \subset \mathcal{M}(k)$.

**Proposition.** *The basic construction gives a bijection*

$$\coprod_{\substack{0 \leq l \leq g-2 \\ \deg L = d-1-(g-2-l)}} Q_{l,L}^*(k) \to \mathcal{M}_1(k),$$

*where the disjoint union is taken over all $l \in [0, g - 2]$ and a set of representatives of all isomorphism classes of line bundles $L$ of degree $d - 1 - (g - 2 - l)$.*

PROOF. By **5**, the map is a surjection. Now suppose that $(l, L, V \subset F_*L)$ and $(l', L', V' \subset F_*L')$ give the same point in $\mathcal{M}_1(k)$, i.e. $V \simeq V'$. Since the unstable bundle $F^*V$ has a unique quotient line bundle of degree $< \deg(V)/2$ (i.e. the second graded piece of the Harder-Narasimhan filtration), which is isomorphic to $L$, we must have $L = L'$. Consider the diagram

$$\begin{array}{ccc} F^*V & \longrightarrow & L \\ \wr \downarrow & & \| \\ F^*V' & \longrightarrow & L', \end{array}$$

where the vertical arrow is induced from an isomorphism $V \xrightarrow{\sim} V'$ and the horizontal arrows are the unique quotient maps. This diagram is commutative up to a multiplicative scalar in $k^*$. By adjunction, $V \hookrightarrow F_*L$ and $V' \hookrightarrow F_*L$ have the same image. In other words, $V = V'$ as sub-modules of $F_*L$. This proves the injectivity of the map. ∎

**7 Moduli space.** To ease the notation, let $d_l = d - 1 - (g - 2 - l)$. Let $\mathrm{Pic}^{d_l} X$ be the moduli space of line bundles of degree $d_l$ on $X$, and $\mathcal{L} \to \mathrm{Pic}^{d_l}(X) \times X$ the universal line bundle.

By [FGA, 3.2], there is a scheme $\mathcal{Q} = \mathcal{Q}_l = \mathrm{Quot}_l\big((\mathrm{id} \times F)_*\mathcal{L}/(\mathrm{Pic}^{d_l}(X) \times X)/\mathrm{Pic}^{d_l} X\big) \xrightarrow{\pi}$ $\mathrm{Pic}^{d_l} X$ such that $\mathcal{Q}_x$ (the fiber at $x$) is $Q_{\mathcal{L}_x}$ for all $x \in (\mathrm{Pic}^{d_l} X)(k)$. By the same argument as before, there is an open sub-scheme $\mathcal{Q}^* \subset \mathcal{Q}$ such that $\mathcal{Q}^*_x = Q^*_{\mathcal{L}_x}$ for all $x \in \mathrm{Pic}^{d_l}(X)(k)$. The scheme $\mathcal{Q}$ is projective over $\mathrm{Pic}^{d_l}(X)$ ([FGA, 3.2]), hence is proper over $k$. By checking the condition of formal smoothness (cf. [L, 8.2.1]), it can be shown that $\mathcal{Q}$ is smooth over $\mathrm{Pic}^{d_l}(X)$, hence is smooth over $k$.

The coarse moduli scheme $\mathcal{M}$ is canonically stratified by Harder-Narasimhan polygons. Concretely, for $j \geq 0$, let $P_j$ be the polygon from $(0,0)$ to $(1, d+j)$ to $(0, 2d)$. Let $\mathcal{M}_0 = \overline{\mathcal{M}}$, and for $j \geq 1$, let $\mathcal{M}_j(k)$ be the subset of $\overline{\mathcal{M}}(k)$ parametrizing those $V$'s such that the Harder-Narasimhan polygons ([S]) of $F^*V$ lie above or are equal to $P_j$. Notice that $\mathcal{M}_1(k)$ agrees with the one defined in **6**.

As mentioned in the introduction, the existence of a universal bundle on $\mathcal{M}$ would imply that each $\mathcal{M}_j(k)$ is Zariski closed by Shatz's theorem [S]. In general, one can show that $\mathcal{M}_j(k)$ is closed by examining the GIT (geometric invariant theory) construction of $\overline{\mathcal{M}}$. This fact also follows from our basic construction:

**Theorem.** *The subset $\mathcal{M}_j(k)$ is Zariski closed in $\overline{\mathcal{M}}(k)$, hence underlies a reduced closed subscheme $\mathcal{M}_j$ of $\overline{\mathcal{M}}$. The scheme $\mathcal{M}_j$ is proper. The Harder-Narasimhan stratum $\mathcal{M}_j \smallsetminus \mathcal{M}_{j+1}$ is non-empty precisely when $0 \leq j \leq g-1$. For $1 \leq j \leq g-1$, write $l = g - 1 - j$. Then there is a canonical morphism*

$$\mathcal{Q}_l \to \overline{\mathcal{M}}$$

*which has scheme-theoretic image $\mathcal{M}_j$ and induces a bijection from $\mathcal{Q}^*_l(k)$ to $\mathcal{M}_j(k) \smallsetminus \mathcal{M}_{j+1}(k)$.*

PROOF. Suppose $0 \leq l \leq g-2$ and $j + l = g - 1$. The universal object $\mathcal{V} \to \mathcal{Q}_l \times X$ is a family of stable vector bundles on $X$. This induces a canonical morphism $\mathcal{Q}_l \to \overline{\mathcal{M}}$. The image of $\mathcal{Q}_l(k)$ is precisely $\mathcal{M}_j(k)$ by (the proof of) Corollary 5.2. Since $\mathcal{Q}_l$ is proper, $\mathcal{M}_j$ is proper and closed in $\overline{\mathcal{M}}$. The rest of the proposition follows from **6** and **5**, and the fact that $\mathcal{Q}^*_l(k)$ is non-empty for $0 \leq l \leq g - 2$ (see Lemma 9.3). ∎

**8 Remark.** By a theorem of Nagata ([LN], [HN]), $s(V) \leq g$ for all $V$. Therefore, $s(V) \leq g$ if $\deg V \equiv g \pmod{2}$, and $s(V) \leq g - 1$ if $\deg V \not\equiv g \pmod{2}$. By Proposition 3, $V = F_*L$ achieves the maximum value of $s$ among rank-2 vector bundles of the same degree.

By the preceding theorem, vector bundles of the form $V = F_*L$ are precisely members of the smallest non-empty Harder-Narasimhan stratum $\mathcal{M}_{g-1}$. Therefore, in a sense $V$ is most stable yet $F^*V$ is most unstable. More generally, for $1 \leq j \leq g - 1$, we have (from **4**)

$$s(\mathcal{M}_j(k)) \geq \begin{cases} j & \text{if } d \equiv j \pmod{2}, \\ j+1 & \text{if } d \not\equiv j \pmod{2}. \end{cases}$$

Therefore, high unstablility of $F^*V$ implies high stability of $V$.

**9  Irreducibility.** We will make use of the following simple lemma.

**9.1  Lemma.** *Let $Y$ be a proper scheme over $k$, $S$ an irreducible scheme of finite type over $k$ of dimension $s$, $r$ an integer $\geq 0$, and $f : Y \to S$ a surjective morphism. Suppose that all fibers of $f$ are irreducible of dimension $r$. Then $Y$ is irreducible of dimension $s + r$.* ∎

**9.2  Lemma.** *The scheme $\mathcal{Q} = \mathcal{Q}_l$ is irreducible of dimension $2l + g$.*

PROOF. There is a surjective morphism ([FGA, §6])

$$\delta : \mathcal{Q} \to \mathrm{Div}^l(X) = \mathrm{Sym}^l(X), \qquad q \mapsto \sum_{P \in X(k)} \mathrm{length}_{\mathcal{O}_P}\big((F_*\mathcal{L}_{\pi(q)})/\mathcal{V}_q\big) \cdot P.$$

The morphism $\mathcal{Q} \to \mathrm{Div}^l(X) \times \mathrm{Pic}^{d_l}(X)$ is again a surjection. The fibers are irreducible schemes of dimension $l$ according to the last lemma of [MX]. Since $\mathcal{Q}$ is proper, the result follows from Lemma 9.1. ∎

**9.3  Lemma.** *$\mathcal{Q}^*$ is open and dense in $\mathcal{Q}$.*

PROOF. By the construction in **6** and **7**, $\mathcal{Q}^*$ is open in $\mathcal{Q}$. Since $\mathcal{Q}$ is irreducible of dimension $2l + g$, it suffices to show that $\mathcal{Q}^*$ is non-empty. We will do more by exhibiting an open subset of $\mathcal{Q}^*$ of dimension $2l + g$.

Indeed, let $B(X, l) \subset \mathrm{Div}^l(X)$ be the open sub-scheme parametrizing multiplicity-free divisors of degree $l$, also known as the configuration space of unordered $l$ points in $X$. Let $U$ be the inverse image of $B(X, l) \times \mathrm{Pic}^{d_l}(X)$ under $\mathcal{Q}^* \to \mathrm{Div}^l(X) \times \mathrm{Pic}^{d_l}(X)$. A quick calculation shows that each fiber of $U \to B(X, l) \times \mathrm{Pic}^{d_l}(X)$ is isomorphic to $\mathbb{A}^l$. Therefore, $U$ is an open subset of $\mathcal{Q}^*$ of dimension $2l + g$. ∎

**9.4  Theorem.** *For $1 \leq j \leq g - 1$, $\mathcal{M}_j$ is proper, irreducible, and of dimension $g + 2(g - 1 - j)$. In particular, $\mathcal{M}_1$ is irreducible and of dimension $3g - 4$.* ∎

**10  Fixing the determinant.** Fix a line bundle $\xi$ of degree $d$. Let $\overline{\mathcal{M}}(\xi) \subset \overline{\mathcal{M}}$ be the closed sub-scheme of $\overline{\mathcal{M}}$ parametrizing those $V$'s with $\det(V) = \xi$. Let $\mathcal{M}_j(\xi) = \overline{\mathcal{M}}(\xi) \cap \mathcal{M}_j$ for $j \geq 0$.

**Remark.** For $1 \leq j \leq g - 1$, $\dim \mathcal{M}_j(\xi) = 2(g - 1 - j)$. In particular, $\dim \mathcal{M}_1(\xi) = 2(g - 2)$.

PROOF. Since $\mathcal{M}_j(\xi)$ is nothing but the fiber of the surjective morphism $\det : \mathcal{M}_j \to \mathrm{Pic}^d(X)$, it has dimension $2(g - 1 - j)$ for a dense open set of $\xi \in \mathrm{Pic}^d(X)(k)$. However, $\mathcal{M}_j(\xi_1)$ is isomorphic to $\mathcal{M}_j(\xi_2)$ for all $\xi_1, \xi_2 \in \mathrm{Pic}^d(X)(k)$, via $V \mapsto V \otimes L$, where $L^2 \simeq \xi_2 \otimes \xi_1^{-1}$. Thus the remark is clear. ∎

A slight variation of the above argument shows that $\mathcal{M}_j(\xi)$ is irreducible. Alternatively, assume $1 \leq j \leq g - 1$. Let $l = g - 1 - j$ and let $\mathcal{Q}(\xi) = \mathcal{Q}_l(\xi)$ be the inverse image of $\xi$ under $\mathcal{Q} \to \text{Pic}^d(X)$, $q \mapsto \det(\mathcal{V}_q)$. Since $\det(\mathcal{V}_q) = B \otimes \mathcal{L}_{\pi(q)} \otimes \mathcal{O}(-\delta(q))$, the morphism $\det : \mathcal{Q} \to \text{Pic}^d(X)$ factors as

$$\mathcal{Q} \to \text{Div}^l(X) \times \text{Pic}^{d_l}(X) \xrightarrow{\psi} \text{Pic}^d(X),$$

where $\psi$ is $(D, L) \mapsto B \otimes L \otimes \mathcal{O}(-D)$. It is clear that $\psi^{-1}(\xi)$ is isomorphic to $\text{Div}^l(X)$, and hence is an irreducible variety.

The fibers of $\mathcal{Q}(\xi) \to \psi^{-1}(\xi)$ are just some fibers of $\mathcal{Q} \to \text{Div}^l(X) \times \text{Pic}^{d_l}(X)$; hence they are irreducible of dimension $l$ as in the proof of Lemma 9.2. Being a closed sub-scheme of $\mathcal{Q}$, $\mathcal{Q}(\xi)$ is proper, thus, irreducible by Lemma 9.1. Now it is easy to deduce

**Theorem.** *The scheme $\overline{\mathcal{M}}(\xi)$ admits a canonical stratification by Harder-Narasimhan polygons*

$$\emptyset = \mathcal{M}_g(\xi) \subset \mathcal{M}_{g-1}(\xi) \subset \cdots \subset \mathcal{M}_0(\xi) = \overline{\mathcal{M}}(\xi),$$

*with $\mathcal{M}_j(\xi)$ non-empty, proper, irreducible, and of dimension $2(g - 1 - j)$ for $1 \leq j \leq g - 1$.* ∎

**11  A variant.** Let $\mathcal{M}'(k)$ be the subset of $\overline{\mathcal{M}}(k)$ consisting of those $V$ such that $F^*V$ is not stable. Clearly, $\mathcal{M}'(k) \supset \mathcal{M}_1(k)$.

By Corollary 5.1, the closed subset $\mathcal{M}^{ns}(k) = \overline{\mathcal{M}}(k) \smallsetminus \mathcal{M}(k)$ is contained in $\mathcal{M}'(k) \smallsetminus \mathcal{M}_1(k)$. On the other hand, if $V \in \mathcal{M}'(k) \smallsetminus \mathcal{M}^{ns}(k)$, the argument of 5 shows that there is a line bundle $L$ of degree $d$ such that $V \hookrightarrow F_*L$ is a sub-module of co-length $\leq g - 1$. Conversely, the argument of 4 shows that if $V$ is of co-length $\leq g - 1$ in $F_*L$ for some $L$ of degree $d$, then $V \in \mathcal{M}'(k)$.

Thus we conclude that $\mathcal{M}'(k)$ is the union of $\mathcal{M}^{ns}(k)$ and the image $\mathcal{M}'_0(k)$ of $\mathcal{Q}_{g-1}(k)$ for a suitable morphism $\mathcal{Q}_{g-1} \to \overline{\mathcal{M}}$, where $\mathcal{Q}_{g-1}$ is defined in 7. It follows that $\mathcal{M}'_0(k)$ and $\mathcal{M}'(k)$ are Zariski closed in $\overline{\mathcal{M}}(k)$, hence are sets of $k$-points of reduced closed sub-scheme $\mathcal{M}'_0$ and $\mathcal{M}'$ of $\overline{\mathcal{M}}$.

**Theorem.** *The scheme $\mathcal{M}'_0$ is irreducible of dimension $3g - 2$. It contains two disjoint closed subsets: $\mathcal{M}'_0 \cap \mathcal{M}^{ns}$, which is irreducible of dimension $2g - 1$ when $d$ is even and empty when $d$ is odd, and $\mathcal{M}_1$, which is irreducible of dimension $3g - 4$.*

**Remark.** $\mathcal{M}' \smallsetminus \mathcal{M}_1$ is the first stratum in the $s$-stratification ([LN]) which is not a Harder-Narasimhan stratum. The other $s$-stratas are more complicated and not pursued here.

PROOF. Since $\mathcal{Q}_{g-1}$ is irreducible, $\mathcal{M}'_0$ is irreducible. We now analyze $\mathcal{M}'_0 \cap \mathcal{M}^{ns}$. Suppose that $V \in \mathcal{M}'_0(k) \cap \mathcal{M}^{ns}(k)$. Then $d = \deg V$ is even and there exists $L$ of degree $d$ such that $V$ is a sub-module of $F_*L$ of co-length $g - 1$. By assumption, there is a sub-bundle $M$ of $V$ of degree $d/2$. Adjunction applied to the composition $M \hookrightarrow V \hookrightarrow F_*L$ provides a non-zero morphism $F^*(M) = M^2 \to L$. This implies that $M^2 \simeq L$. We may assume that $L = M^2$. Since there is only one (modulo $k^*$) non-zero morphism $M^2 \to L$, there is only one non-zero morphism $M \to F_*(F^*M)$. By [JX, §2], this morphism is part of an exact sequence $0 \to M \to F_*(F^*M) \to M \otimes B \to 0$. Thus to have $V$ is to have a sub-module of $M \otimes B$ of co-length $g - 1$. Conversely, starting with a sub-module of $M \otimes B$ of co-length $g - 1$, we obtain a vector bundle $V \in \mathcal{M}_0'(k) \cap \mathcal{M}^{ns}(k)$ as the inverse image of that sub-module in $F_*(F^*M)$.

The sub-modules of $M \otimes B$ of co-length $g - 1$ are of the form $M \otimes B \otimes \mathcal{O}(-D)$ for $D \in \mathrm{Div}^{g-1}(X)(k)$. Thus there is a morphism $\pi' : \mathcal{Q}' = \mathrm{Div}^{g-1}(X) \times \mathrm{Pic}^{d/2}(X) \to \overline{\mathcal{M}}$ inducing a surjection $\mathcal{Q}'(k) \to \mathcal{M}_0'(k) \cap \mathcal{M}^{ns}(k)$. We claim that this morphism is generically finite of separable degree at most 2. This claim implies that $\mathcal{M}_0' \cap \mathcal{M}^{ns}$ is irreducible of dimension $2g - 1$.

Indeed, there is an open subset $U$ of $\mathrm{Div}^{g-1}(X)(k)$ such that if $D, D' \in U$ are distinct, then $D \not\sim D'$. We now show that $\pi'|(U \times \mathrm{Pic}^{d/2}(X)(k))$ is at most 2-to-1. Suppose that $D \in U$, $M \in \mathrm{Pic}^{d/2}(X)(k)$, and $\pi'(D, M) = V$. Then $V$ has at most two isomorphism classes of rank-1 sub-bundles of degree $d/2$, and $M$ is one of them. After obtaining $M$, one can determine $D$ uniquely by the condition $\det(V) \simeq M^2 \otimes B \otimes \mathcal{O}(-D)$. This proves the claim.

Next, we consider the morphism $\mathcal{Q}_{g-1} \to \mathcal{M}_0'$. It induces a surjection $\mathcal{Q}_{g-1}^*(k) \twoheadrightarrow \mathcal{M}_0'(k) \smallsetminus \mathcal{M}_1(k)$. Again the claim is that the morphism is generically finite of separable degree at most 2. This claim implies that $\mathcal{M}_0'$ is irreducible of dimension $3g - 2$.

Indeed, let $U$ be the open subset of $\mathcal{Q}_{g-1}^*(k)$ consistng of those $q$'s such that $\mathcal{O}(2\delta(q)) \not\simeq \Omega_{X/k}$. Now assume that $q \in U$ gives rise to $V \in \mathcal{M}_0'(k)$. Then there is an exact sequence $0 \to L \otimes B^2 \otimes \mathcal{O}(-2\delta(q)) \to F^*V \to L \to 0$, where $L = \mathcal{L}_{\pi(q)}$. The assumption on $q$ implies that $F^*V$ has at most 2 quotient line bundles of degree $d$, say $F^*V \to L_1$ and $F^*V \to L_2$. Then $q$ must be one of the two data $V \hookrightarrow F_*L_1$ or $V \hookrightarrow F_*L_2$ provided by adjunction. This proves the claim. ∎

**12 Example.** When $g = 2$, $\mathcal{M}_1(\xi)$ is a single point, corresponding to the vector bundle $F_*(\xi \otimes B^{-1})$.

When $\xi = B$, this refines a result of Joshi and one of us [JX, 1.1], which says that $\mathcal{M}_1(\xi)$ is a single $\mathrm{Pic}(X)[2]$-orbit.

When $\xi = \mathcal{O}_X$, this extends a theorem of Mehta [JX, 3.2], which states that there are only finitely many rank-2 semi-stable vector bundles $V$'s on $X$ with $\det(V) = \mathcal{O}_X$ and $F^*V$ not semi-stable when $p \geq 3$, $g = 2$. We now have this result for $p = 2, g = 2$ with the stronger conclusion of uniqueness.

**13 Erratum for [JX].** We correct a minor error in the statement of [JX, Theorem 1.1]. The expression "$V_1 \in \mathrm{Ext}^1(L_\theta, \mathcal{O}_X)$" should be replaced by "$V_1 \in S_\theta$" (the original version is valid when $L_\theta = B$). Also, $\Omega$ should be replaced by $L_\theta$.

**14 References.**

CURRENT ADDRESSES:

Jiu-Kang Yu
Department of Mathematics
University of Maryland
College Park, MD 20740
U.S.A.
*Email:* `yu@math.umd.edu`

Eugene Z. Xia
Department of Mathematics & Statistics
University of Massachusettes
Amherst, MA 01003
U.S.A.
*Email:* `xia@math.umass.edu`